\documentstyle[amscd,amssymb,epsf,12pt]{amsart}
\newtheorem{thm}{Theorem}[section]
\newtheorem{lem}[thm]{Lemma}

\newtheorem{prop}[thm]{Proposition}
\newtheorem{cor}[thm]{Corollary}

\theoremstyle{definition}

\newtheorem{rem}[thm]{Remark}
\newtheorem{nota}[thm]{Notation}

\makeatletter
\def\theequation{\@arabic\c@equation}
\makeatother
\numberwithin{equation}{thm}

\def\eqq{\stackrel{\sim}{=}}

\def\Q{{\Bbb Q}}
\def\Z{{\Bbb Z}}

\def\P{{\Bbb P}}
\def\C{{\Bbb C}}

\def\C{{\Bbb C}}

\def\P{{\Bbb P}}

\def\Z{{\Bbb Z}}
\def\Q{{\Bbb Q}}

\def\delt{\bigtriangleup}

\def\({\left(}
\def\){\right)}

\def\part{P(n)}

\def\<{\langle}
\def\>{\rangle}

\def\scup{\mathbin{\text{\scriptsize$\cup$}}}
\def\scap{\mathbin{\text{\scriptsize$\cap$}}}

\begin{document}

\title[The quantum cohomology of blow-ups of $\P^2$]
{The quantum cohomology of blow-ups of $\P^2$ and enumerative geometry}
\author{L. G\"ottsche }
\author{R. Pandharipande}
\date{11 November 1996}
\thanks{The second author was partially supported by an NSF post-doctoral
fellowship.}
\address{Mittag-Leffler-Institute\\ Aurav\"agen 17\\
Djursholm, Sweden}
\email{gottsche@@ml.kva.se,
 pandhari@@ml.kva.se}

 \maketitle\



\section{Introduction}
\label{intro}
The enumerative geometry of curves in algebraic varieties has taken a new
direction with the appearance of Gromov-Witten 
invariants and quantum cohomology.
Gromov-Witten invariants originate in symplectic geometry and were
first defined in terms of pseudo-holomorphic curves.
In algebraic geometry,  these invariants are defined using moduli spaces 
of stable maps.

Let $X$ be a nonsingular
projective variety   over $\C$. Let $\beta\in H_2(X,\Z)$.
In  [K-M], the moduli space $\overline M_{0,n}(X,\beta)$ of stable
$n$-pointed genus $0$ maps is defined. This moduli space parametrizes
the data 
$[\mu:C\rightarrow X,p_1,\dots,p_n]$ where $C$ is a connected, reduced,
(at worst) nodal curve of genus 0,
$p_1,\ldots,p_n$ are nonsingular points of $C$, and $\mu$ is a morphism.
$\overline M_{0,n}(X,\beta)$
is equipped with
$n$ morphisms $\rho_1,\ldots,\rho_n$ to $X$ where
$\rho_i ([\mu:C\rightarrow X,p_1,\dots,p_n]) =\mu(p_i)$.

$X$ is a convex variety if $H^1(\P^1,f^*(T_X))=0$ 
for all maps $f:\P^1 \rightarrow X$. In this case,  
$\overline M_{0,n}(X,\beta)$ is a  projective scheme 
  of pure expected dimension equal to  
$$dim(X)+n-3+\int_\beta c_1(T_X)$$
with only finite quotient singularities.
Given classes $\gamma_1,\dots,\gamma_n$ in $H^*(X,\Z)$, the 
Gromov-Witten invariants $I_\beta(\gamma_1\dots\gamma_n)$ are defined
by:
$$
I_\beta(\gamma_1\dots\gamma_n)=
\int_{\overline M_{0,n}(X,\beta)}
\rho_1^*(\gamma_1)\scup \ldots \scup \rho_n^*(\gamma_n).
$$
The intuition behind these invariants is as follows.
If the $\gamma_i$ are the cohomology classes of subvarieties
$Y_i\subset X$ in general position,  then 
$I_\beta(\gamma_1\dots\gamma_n)$ should count the (possibly virtual)
number of irreducible
rational  curves $C$ in $X$ of homology class $\beta$ which intersect
all the $Y_i$. In case $X$ is a homogeneous space, a
correspondence between the Gromov-Witten invariants and
the enumerative geometry of rational curves in $X$ can be proven
by transversality arguments (see [F-P]).  

One can use the Gromov-Witten invariants to define 
the big quantum cohomology ring $QH^*(X)$ of $X$. 
The associativity of this ring yields relations 
among the invariants $I_\beta(\gamma_1\dots\gamma_n)$ which often 
are sufficient to determine
them all recursively from a few basic ones.   
The model case for this approach
is the recursive determination of the numbers $N_d$ of nodal rational curves
of degree $d$ in the projective plane [K-M], [R-T].

If $X$ is not convex, the moduli space $\overline M_{0,n}(X,\beta)$
will not in general have the expected dimension.
Recently, Gromov-Witten invariants have been defined
and proven to satisfy basic geometric properties
via the construction of  virtual fundamental classes of 
the expected dimension [B-F],  [B],  [L-T 2]
and, in the symplectic context,  [L-T 3], [F-O], [S].
In particular, these Gromov-Witten invariants
have been proven to satisfy the axioms of
[K-M], [B-M]. Therefore, they again define an associative
quantum cohomology ring $QH^*(X)$.

The aim of this paper is to study the Gromov-Witten invariants 
of the blow-up $X_r$ 
of $\P^2$ in a finite set $x_1,\ldots x_r$ of points and
to give enumerative applications. $X_r$ is a particularly simple example of
a nonconvex variety, so this study (at least in the context of algebraic
geometry) neccessitates the use of the above constructions.
Let $S$ be a nonsingular, rational, projective surface. 
$S$ is either deformation equivalent to $\P^1 \times
\P^1$
or to $X_{r(S)}$ where $r(S)+1= rank( A^1(S))$.
Together with the invariants of $\P^1 \times \P^1$, the
Gromov-Witten invariants of $X_r$ therefore
determine the invariants
of all these rational surfaces (the invariants are
constant in flat families of nonsingular varieties).
For enumerative applications,
it is necessary to consider the blow-up $X_r$ of
$\P^2$ in a finite set of {\em general} points. 

Let $H$ be the pull-back to $X_r$
of the hyperplane class in $\P^2$, and let
$E_1,\ldots,E_r$  be the exceptional divisors. Our aim is to count
the number of irreducible rational curves $C$ in $X_r$ of class
$dH-\sum_{i=1}^r a_iE_i$ passing through
$3d-\sum_{i=1}^r a_i -1$ general points. By associating to a curve 
in $\P^2$ its strict transform in $X_r$, this number
can also be interpreted as the number of irreducible rational curves
in $\P^2$ having singularities of 
order $a_i$ at the (fixed) general points $x_i$
and passing  through $3d-\sum_{i=1}^r a_i -1$ other general points.

The paper is naturally divided  into two parts. First,
we use the associativity
of the quantum product to show that the Gromov-Witten invariants of
$X_r$ can
be computed from simple initial values by means of explicit recursion
relations. There are $r+1$ initial values required
for $X_r$:
\begin{enumerate}
\item[(i)] The number of lines in the plane passing through 2
points, $N_{1,(0, \ldots, 0)}=1$.
\item[(ii)] The number
 of curves in the exceptional class $E_i$,
 $N_{0, -[i]}=1$.
\end{enumerate}
The relations are then used to prove  properties 
of these invariants. 

In the second half of the paper, the  
enumerative significance of the invariants
is investigated. Our main tool is a degeneration
argument in which the points $x_i$ are specialized to lie on a 
nonsingular cubic in $\P^2$.
The idea of using such degenerations is due independently to
J. Koll\'ar and, in joint work, to    
L. Caporaso and J. Harris [C-H]. 
For a general blow-up $X_r$, the Gromov-Witten
invariants are proven to be a count (with possible multiplicities)
of the finite number of solutions to the corresponding
enumerative problem on $X_r$. 
Let $\beta=dH-\sum_{i=1}^{r} a_iE_i$ be a class in $H_2(X_r, \Z)$.
If the expected  dimension
of the moduli space 
$\overline M_{0,0}(X_r,\beta)$ is strictly positive or
if there exists a multiplicity $a_i \in \{1,2\}$, then
the corresponding Gromov-Witten invariant is proven to
be an actual count of the number of irreducible,
degree $d$, rational plane curves of multiplicity $a_i$
at the (fixed) general points $x_i$ which pass through
$3d-\sum_{i=1}^{r} a_i -1$ other general points.
In the Del Pezzo case ($r\leq 8$), all
invariants are shown to be enumerative (see also [R-T]).
A basic symmetry of the
Gromov-Witten invariants of the spaces $X_r$ obtained from the 
classical
Cremona transformation is discussed in section 
\ref{cremmy}. 
These considerations show that for $d\leq 10$,
the Gromov-Witten invariants always
coincide with enumerative geometry.  
Tables 
of these invariants in low degrees are given in
section \ref{tbls}.

In [K-M], an associativity equation for
Del Pezzo surfaces (corresponding to our
relation $R(m)$) is derived.
The small quantum cohomology ring of Del Pezzo surfaces
is studied in [C-M].
In section 11 of [C-M], the associativity of the
small quantum product on $X_r$ is used to
derive some relations among the Gromov-Witten
invariants of these surfaces.   
The invariants of $\P^2$ blown-up 
in a point are computed in [C-H2], [G], and [K-P].
In [G], A. Gathmann 
computes more generally the invariants
of the 
blow-up of $\P^n$ in a point and studies their
enumerative significance.
In [D-I], the Gromov-Witten invariants
of $X_6$ are computed via associativity. Our recursive
strategy for $X_6$ differs.

The first author would like to thank K. Hulek who suggested
to him the possibility of studying the quantum cohomology of blow-ups.
He also thanks the University of Chicago and in particular 
W. Fulton for an invitation which made this collaboration possible.

The second author would like to thank J. Koll\'ar
for many remarkable insights on these questions. In particular, the
main immersion result (Lemma \ref{lemmy}) was first observed by
him. Thanks are due to J. Harris for conversations in which 
the degeneration argument was explained.

Both authors thank the Mittag-Leffler-Institut for support.

\section{Notation and Background material}
\label{notabackg}
Let $X$ be a  nonsingular projective variety. 
Assume for simplicity that the Chow and homology rings of
$X$ coincide. Let $dim(X)$ be the complex dimension.
Denote by 
$\alpha\scup \beta$ the cup product of classes $\alpha,\beta\in H^*(X,\Z)$
and let $(\alpha\cdot \beta)=\int_X\alpha\scup\beta$. 
By definition, $(\alpha\cdot \beta)$ is zero if $\alpha\in H^{2i}(X,\Z)$, 
$\beta\in H^{2j}(X,\Z)$, and $i+j\ne dim(X)$.

We  recall the definition of quantum cohomology from
 [K-M] in a slightly modified form for
nonconvex varieties. Let $B\subset H_2(X,\Z)$ be the semigroup 
of non-negative linear combinations of classes of algebraic curves. 
Let $\beta\in H_2(X, \Z)$. Let $n_\beta=dim(X)+\int_{\beta} 
c_1(T_X) -3$.
Let $n\geq 0$.
For classes   
$\gamma_i\in H^{2j_i}(X,\Z)$ 
 with $\sum_{i=1}^{n} j_i=n_{\beta}+n$,
let $I_\beta(\gamma_1\dots \gamma_n)$
be the corresponding Gromov-Witten invariant:
$$
I_\beta(\gamma_1\dots\gamma_n)=
\int_{[\overline M_{0,n}(X,\beta)]}
\rho_1^*(\gamma_1)\scup \ldots \scup \rho_n^*(\gamma_n)
$$
where $[\overline M_{0,n}(X,\beta)]$
is the virtual fundamental class. Note
if $n_\beta=0$ and $n=0$, then  $I_\beta$ is just the degree
of the fundamental class.
Kontsevich and Manin introduced a set of axioms 
for the Gromov-Witten invariants
which have now been established for nonsingular varieties
(see section \ref{intro}).
If $\overline M_{0,n}(X,\beta)$ is empty, then
 $I_\beta(\gamma_1\dots \gamma_n)=0$. In particular,
all invariants vanish for $\beta\not\in B$.
Let $T_0=1,T_1,\ldots,T_m$ be a homogeneous $\Z$-basis for $H^*(X,\Z)$.
We assume that $T_1,\ldots,T_p$ form a basis of $H^2(X,\Z)=Pic(X)$.
We denote by $T_i^\vee$ the corresponding elements of the dual basis:
$(T_i\cdot T_j^\vee)=\delta_{ij}$.
Denote by $(g_{ij})$ the matrix of intersection numbers 
$(T_i\cdot T_j)$ and by $(g^{ij})$ the inverse matrix.
For variables $y_0,q_1,\ldots,q_p,y_{p+1},\ldots,y_{m}$
(which we also abbreviate as $q,y$), define 
the formal power series 
\begin{equation}
\label{potfun}
\Gamma(q,y)=\sum_{n_{p+1}+\ldots+n_{m}\ge 0}\sum_{\beta\in B \setminus
\{0\}}
I_\beta(T_{p+1}^{n_{p+1}}\cdots T_m^{n_m})
q_1^{\int_\beta  T_1}\cdots 
q_p^{\int_\beta  T_p}
\frac{y_{p+1}^{n_{p+1}}\cdots y_{m}^{n_m}}
{n_{p+1}!\cdots n_{m}!}
\end{equation}
in the ring
$$\Q[[q, q^{-1}, y]]=\Q[[y_0, q_1, \ldots, q_p, q_1^{-1},
\ldots, q_p^{-1}, y_{p+1}, \ldots, y_m]].$$

In case $X$ is a homogeneous space,
the substitution $q_i=e^{y_i}$ in (\ref{potfun}) yields
a formal power series which equals 
the quantum part of the potential function of [K-M] modulo
a quadratic polynomial in the variables $y_1, \ldots, y_m$.
The form (\ref{potfun}) of the potential function
is chosen to avoid convergence issues
in the nonconvex case.
Let
$$
\partial_i=
\begin{cases}
q_i\frac{\partial}{\partial q_i}& i=1,\ldots,p\\
\frac{\partial}{\partial y_i} & i=0,p+1,\ldots,m
\end{cases}
$$
and denote 
$f_{ijk}=\partial_i\partial_j\partial_k f$ for $f\in \Q[[q,q^{-1},y]]$.
Define a $\Q[[q,q^{-1},y]]$-algebra structure on the
free $\Q[[q,q^{-1},y]]$-module generated by $T_0,\ldots,T_m$ by:
$$T_i*T_j=T_i\scup T_j+\sum_{e,f=0}^m \Gamma_{ije} g^{ef}T_f.$$
By definition, this is the quantum cohomology ring
of $X$, $QH^*(X)$.

We sketch the proof of the associativity of this quantum product
following [K-M] and  [F-P]. 
First, a formal calculation (using the axiom of
divisor) yields:
\begin{equation}
\label{partpotfun}
\Gamma_{ijk}=\sum_{n \ge 0}\sum_{\beta\in B \setminus
\{0\}} \frac{1}{n!}
I_\beta( \gamma^n \cdot T_i T_j T_k)
q_1^{\int_\beta  T_1}\cdots 
q_p^{\int_\beta  T_p},
\end{equation}
where
$\gamma= y_{p+1} T_{p+1} +\ldots +y_m T_m$
and the $\Q[[y_0,y_{p+1},\ldots,y_{m}]]$-linear extension of 
$I_\beta$ is used.
Define the symbol $\Phi_{ijk}$ by $\Phi_{ijk}= I_0(T_i T_j 
T_k)+ \Gamma_{ijk}$. In case $X$ is homogeneous, $\Phi_{ijk}$
is the partial derivative of the full potential function.
The $*$-product can be expressed by:
$T_i * T_j= \sum_{e,f=0}^m \Phi_{ije} g^{ef} T_f.$
Let
$$F(i,j|k,l)=\sum_{e,f=0}^{m}\Phi_{ije}g^{ef}\Phi_{fkl}.$$
Associativity is now equivalent to  $F(i,j|k,l)=F(j,k|i,l)$.
Following [F-P], we let 
\begin{equation}
\label{ggdef}
G(i,j|k,l)_{\beta,n}=\sum {n\choose n_1}g^{ef}
I_{\beta_1}(\gamma^{n_1}\cdot T_i T_j T_e)
I_{\beta_2}(\gamma^{n_2}\cdot T_k T_l T_f)
\end{equation}
where the sum runs over all $n_1,n_2\geq 0$ with $n_1+n_2=n$ and 
all $\beta_1,\beta_2\in B$ with $\beta_1+\beta_2=\beta$.
As before,
 $\gamma=y_{p+1} T_{p+1} +\ldots +y_m T_m$.
A calculation using equations (\ref{partpotfun}) and (\ref{ggdef})
yields:
$$
F(i,j|k,l)=\sum_{\beta\in B} q_1^{\int_\beta T_1}\dots
q_{p}^{\int_\beta  T_p} \sum_{n\ge 0}\frac{1}{n!} G(i,j|k,l)_{\beta,n}.
$$
On the other hand, we can use the splitting axiom and linear equivalence on
$\overline M_{0,4}=\P^1$ to see 
that $G(i,j|k,l)_{\beta,n}=G(j,k|i,l)_{\beta,n}$,
and thus the associativity follows.

\section{Quantum cohomology of blow-ups of $\P^2$}
\label{qcb}
\begin{nota}
Let $r\geq 0$. Let $X_r$ be the blowup of $\P^2$ in $r$ general  points
$x_1,\ldots,x_r$.
Denote by $H\in H_2(X,\Z)$  the hyperplane class and by $E_i$,  for 
$i=1,\ldots ,r$,
the exceptional divisors.
Let  $m=r+2$. Let $T_0=1$. Let 
$T_1$, $T_{i+1}$ ( for $i=1,\ldots,r$), and  $T_m$  
be the Poincar\'e dual cohomology classes of 
$H$, $E_i$ and the class of a point respectively. 
Let  $\epsilon_1=1$ and $\epsilon_i=-1$ for 
$i=2,\ldots ,r+1$. Then, $T_0^\vee=T_m$ and $T_i^\vee=\epsilon_iT_i$
for $i=1,\ldots ,r+1$.
 For an $r$-tuple $\alpha=(a_1,\ldots,a_r)$ of integers, denote by
$(d,\alpha)$ the class
$dH-\sum_{i=1}^r a_i E_i$. Let     $|\alpha|=\sum_i a_i$, and let
 $n_{d,\alpha}=3d-|\alpha|-1$ be the expected dimension of the
moduli space $\overline M_{0,0}(X_r,(d,\alpha))$. 
If $n_{d,\alpha}\geq 0$,
let
$$
N_{d,\alpha}=I_{(d,\alpha)}(T_m^{n_{d,\alpha}})
$$
be the corresponding Gromov-Witten invariant.
When writing $N_{d,\alpha}$ for $\alpha$ a sequence 
of length $r$, we will always mean the Gromov-Witten invariant  
on $X_r$.

The   components of the finite sequences
$\alpha$, $\beta$, $\gamma$ are   denoted by the corresponding
roman letters $a_i$, $b_i$, $c_i$. 
For any $r$, we write $[i]_r$ for the sequence $(j_1,\ldots,j_r)$
with $j_k=\delta_{ik}$. We just write $[i]$ if $r$ is
understood.
For a sequence $\beta=(b_1,\ldots,b_{r-1})$,
we denote by $(\beta,k)$ the sequence obtained by adding  $b_r=k$.
For a permutation $\sigma$ of $\{1,\ldots,r\}$, denote by 
$\alpha_\sigma$ the sequence $(a_{\sigma(1)},\ldots,a_{\sigma(r)})$.
For an integer $k$,
 we write $\alpha\ge k$ to mean that $a_i\ge k$ for all $i$.
\end{nota}

The invariants $N_{1,(0, \ldots, 0)}$ and 
$N_{0,-[i]_r}$ 
are first determined. A result relating virtual and
actual fundmental classes is needed. Let 
$\overline{M}_{0,0}^*(X,\beta)$ denote the open
locus of automorphism-free maps ($\overline{M}_{0,0}^*(X,\beta)$
is a fine moduli space).
\begin{prop}
\label{fb} 
If $\overline{M}_{0,0}(X, \beta) = \overline{M}_{0,0}^*(X,\beta)$
and the moduli space
is of pure expected dimension, then the virtual fundamental
class is the ordinary scheme theoretic 
fundamental class $[\overline{M}_{0,0}(X, \beta)]$.
\end{prop}
\noindent If, in addition, the expected dimension is 0, 
then the Gromov-Witten invariant $N_{\beta}$ equals
the (scheme-theoretic) length of $\overline{M}_{0,0}(X, \beta)$.
This result is a direct consequence of the construction in [B-F].
\begin{lem} $N_{1,(0, \ldots, 0)}=1$ and $N_{0,-[i]_r}=1$.
\end{lem}
\begin{pf}
A simple check shows that
$\overline  M_{0,2}(X_r,H)=\overline M_{0,2}^*(X_r,H)$. Also,
the moduli space is irreducible of dimension 4 and
(at least) generically nonsingular. 
For
two general points $p_1, p_2 \in X_r$,
$\rho_1^{-1}(p_1)\scap \rho_2^{-2}(p_2)$ 
consists of one reduced point corresponding to preimage
of the  unique line  connecting  the images of $p_1$ and $p_2$ in $\P^2$.
Hence, $N_{1, (0,\ldots,0)}= 1$ by Proposition \ref{fb}.

The moduli space
$\overline M_{0,0}(X_r,(0,-[i])$ 
consists of one automorphism-free map  $\mu: \P^1 \stackrel {\sim}
{\rightarrow} E_i \subset X_r$.
The Zariski tangent space to $\overline M_{0,0}(X_r,(0,-[i])$ 
at $[\mu]$ is $H^0(\P^1, N_{X_r})=0$ where $N_{X_r}\eqq 
{\cal O}_{\P^1}(-1)$
is the normal bundle of the map $\mu$. Hence,
$\overline M_{0,0}(X_r,(0,-[i])$ is nonsingular and
$N_{0,-[i]}=1$ by Proposition \ref{fb}.
\end{pf}

The invariants $N_{d,\alpha}$ will be determined by explicit
recursions. In addition, these Gromov-Witten invariants
will be shown to satisfy the following geometric
properties.
\begin{enumerate}
\item[(P1)]
$N_{0,\alpha}=0$ unless $\alpha=-[i]$ for some $i$.
\item[(P2)]
$N_{d,\alpha}=0$ if $d>0$ and any of the $a_i$ is negative.
\item[(P3)]
$N_{d,\alpha}=N_{d,\alpha_{\sigma}}$ for any permutation $\sigma$ of
$\{1,\ldots,r\}$.
\item[(P4)]
  $N_{d,\alpha}=N_{d,(\alpha,0)}$.
In particular  $N_{d,(0,\ldots,0)}$ 
is the number of rational curves on $\P^2$ passing through $3d-1$ 
general points computed by recursion in [K-M].
\item[(P5)]
If $n_{d,\alpha}>0$, then $N_{d,\alpha}=N_{d,(\alpha,1)}$.
\end{enumerate}

\begin{rem}\label{p11}
Let $Y$ be the blow-up of $\P^1\times\P^1$ in a point
with exceptional divisor $E$, and let $F,G$ be the pullbacks of the classes
of the fibres of the two projections to $\P^1$. 
There is an isomorphism $\phi:X_2\to Y$ with
$\phi_*(H)=F+G-E,$ $\phi_*(E_1)=F-E$, $\phi_*(E_2)=G-E$.
Let $(d,\alpha)$ be given with $r\ge 2$.
If $d-a_1-a_2\ge 0$, then   pushing down   
first to $X_2$ and then further to $\P^1\times \P^1$
gives a bijection between the  irreducible 
rational curves in $|(d,\alpha)|$ on  $X_r$ passing 
through $n_{d,\alpha}$ general points and the 
irreducible rational curves
of bidegree $(d-a_1,d-a_2)$ on $\P^1\times\P^1$, 
with points of multiplicities
$d-a_1-a_2, a_3, \ldots ,a_r$ at $r-1$ general points and passing through 
$n_{d,\alpha}$ other general points.
\end{rem}

We obtain   recursion formulas determining the $N_{d,\alpha}$
from the associativity of the quantum product.
All effective classes $(d,\alpha)$ on $X_r$  satisfy 
$\alpha\leq d$. Therefore, we can write  
$$\Gamma(q,y)=
\sum_{(d,\alpha)} 
N_{d,\alpha} q_1^{d}q_2^{a_1}\dots q_{r+1}^{a_r}
\frac{y_m^{n_{d,\alpha}}}{n_{d,\alpha}!},
$$
where the sum runs over 
all $(d,\alpha)\ne 0$ satisfying $n_{d,\alpha}\geq0$, 
$d\ge 0$, and $\alpha\leq d$.
Let
$\Gamma_{ijk}=\partial_i\partial_j\partial_k\Gamma$ (following the
notation of
section \ref{notabackg}) .
The quantum product of  $T_i$ and $T_j$ is  given by 
$$
T_i*T_j=(T_i\cdot T_j)T_m+\sum_{k=1}^{r+1}\epsilon_k \Gamma_{ijk}T_k+ 
\Gamma_{ijm} T_0.
$$

\begin{lem}
For $i,j,k,l\in \{1,\ldots,m\}$, there is a relation:  
 \begin{align*}\tag{$R_{i,j,k,l}$}
(T_i\cdot T_j)&\Gamma_{klm}-(T_k\cdot T_j)\Gamma_{ilm}+
(T_k\cdot T_l)\Gamma_{ijm}-(T_i\cdot T_l)\Gamma_{kjm}=\\
&=\sum_{s=1}^{m-1} 
\epsilon_s(\Gamma_{jks}\Gamma_{isl}-\Gamma_{ijs}\Gamma_{ksl}).
\end{align*}
\end{lem}
\begin{pf}
We write 
$$
(T_i*T_j)*T_k-(T_k*T_j)*T_i=
\sum_{l=0}^m r_{i,j,k,l} T_l^\vee.
$$
By associativity, we obtain  the relation   $r_{i,j,k,l}=0$.
We show this relation is equivalent to  $(R_{i,j,k,l})$.
We   compute  directly 
\begin{align*}
(T_i&*T_j)*T_k=
(T_i\cdot T_j)T_m*T_k+\sum_{s=1}^{m-1}\epsilon_s\Gamma_{ijs}T_s*T_k
+\Gamma_{ijm}T_k\\
&=\sum_{l=1}^m(T_i\cdot T_j)\Gamma_{klm}T_l^\vee+
\sum_{s=1}^{m-1}\Bigg(\epsilon_s\Gamma_{ijs} 
(T_s\cdot T_k)T_m+\sum_{l=1}^m 
\epsilon_s\Gamma_{ijs}\Gamma_{ksl}T_l^\vee\Bigg)+\Gamma_{ijm}T_k.
\end{align*}
It is easy to see that
\begin{align*}
\Gamma_{ijm}T_k
&= \sum_{l=1}^{m}\Gamma_{ijm} (T_k\cdot T_l)T_l^\vee
+\Gamma_{ijm}\delta_{km} T_0^\vee,\\
\sum_{s=1}^{m-1}\epsilon_s\Gamma_{ijs}(T_s\cdot T_k) T_m&=
\Gamma_{ijk}(1-\delta_{km})T_0^\vee.
\end{align*}
Therefore, the sum of these two terms is just
$ \sum_{l=1}^{m}\Gamma_{ijm} 
\big((T_k\cdot T_l)T_l^\vee+\Gamma_{ijk}T_0^\vee.$
Thus
$$(T_i*T_j)*T_k=\sum_{l=1}^m
\Bigg((T_i\cdot T_j)\Gamma_{klm}+(T_k\cdot T_l)\Gamma_{ijm}
+\sum_{s=1}^{m-1}\epsilon_s\Gamma_{ijs}\Gamma_{ksl}\Bigg)T_l^\vee
+\Gamma_{ijk}T_0^\vee,$$
and the result follows by exchanging the role of $i$ and $k$ and subtracting.
\end{pf}

For the recursive determination of the $N_{d,\alpha}$,  only
the following relations are needed:
\begin{equation}
\tag{$R_{1,1,m,m}$}
\Gamma_{mmm}=\sum_{s=1}^{m-1} \epsilon_s(\Gamma_{1sm}^2-\Gamma_{11s}
\Gamma_{smm}),
\end{equation}
and for all $i=2\ldots r+1$ 
\begin{equation}\tag{$R_{1,1,i,i}$}
\Gamma_{iim}-\Gamma_{11m}=
\sum_{s=1}^{m-1} \epsilon_s(\Gamma_{1is}^2-\Gamma_{11s}
\Gamma_{iis}).
\end{equation}
Note that in case $r=0$, only the relation 
$(R_{1,1,m,m})$ occurs and coincides with that of [K-M].
In the summations below, the following notation is used.
Let the symbol $\vdash\!\! (d,\alpha)$ denote the set of
pairs $\left( (d_1, \beta), (d_2, \gamma) \right)$
satisfying:
\begin{enumerate}
\item[(i)] $(d_1, \beta), (d_2, \gamma)\neq 0$,
\item[(ii)]  $(d_1, \beta)+(d_2,\gamma)=(d, \alpha)$,
\item[(iii)] $n_{d_1,\beta}, n_{d_2, \gamma}\geq 0$, $d_1, d_2 \geq 0,$
$\beta\leq d_1$, and $\gamma\leq d_2$.
\end{enumerate}
The notation $\vdash\!\!(d,\alpha), d_i>0$ will
be used to denote the subset of $\vdash\!\!(d,\alpha)$
satisfying $d_1, d_2>0$.
The binomial coefficient $\binom{p}{q}$ is defined to be
zero if $q<0$ or $p<q$.

\begin{thm}
\label{recurr}
The $N_{d,\alpha}$ are determined  by the initial values:
\begin{enumerate}
\item[(i)] $N_{1,(\underbrace{0,\ldots,0}_r)}=1$, for all $r$,
\item[(ii)] $N_{0,-[i]_r}=1$, for $i\in \{1,\ldots,r\}$,
\end{enumerate}
and the following  recursion relations.
\vspace{+15pt}

\noindent 
If $n_{d,\alpha}\ge 3$, then relation $R(m)$ holds: 
\begin{align*}
 N_{d,\alpha}=
\sum_{\vdash(d,\alpha), d_i>0}&N_{d_1,\beta}N_{d_2,\gamma}
\Big(d_1d_2
-\sum_{k=1}^{r}b_kc_k\Big)
\left(
d_1d_2\binom{n_{d,\alpha}-3}{n_{d_1,\beta}-1}
-d_1^2\binom{n_{d,\alpha}-3}{n_{d_1,\beta}}
\right).
\end{align*}
\vspace{+5pt} 

\noindent
If $n_{d,\alpha}\ge 0$, then 
  for any $i\in \{1,\ldots,r\}$ relation $R(i)$ holds:
\begin{align*}
d^2a_iN_{d,\alpha}&=
(d^2-(a_i-1)^2)N_{d,\alpha-[i]} \\
&+\sum_{\vdash(d,\alpha-[i]), d_i>0}
N_{d_1,\beta}N_{d_2,\gamma}
\Big(d_1d_2
-\sum_{k=1}^{r}b_kc_k\Big)
\left(d_1d_2b_ic_i-
d_1^2c_i^2\right)
\binom{n_{d,\alpha}}{n_{d_1,\beta}}.
\end{align*}

\noindent
Furthermore, the properties (P1)-(P5) hold.
\end{thm}
\begin{pf}
From the relation $(R_{1,1,i+1,i+1})$ above, we get immediately
(for $n_{d,\alpha}\ge 1$) the 
recursion formula $R(i)^*$:
\begin{align*}
(a_i^2-d^2)N_{d,\alpha}=
\sum_{\vdash(d,\alpha)} N_{d_1,\beta}N_{d_2,\gamma}
\Big(d_1d_2
-\sum_{k=1}^{r}b_kc_k\Big)
\left(d_1d_2b_ic_i-
d_1^2c_i^2\right)
\binom{n_{d,\alpha}-1}{n_{d_1,\beta}}.
\end{align*}

We now show  property (P1).
If  $N_{0,\alpha}\ne 0$, then $(0,\alpha)$ is effective
and therefore $\alpha\leq 0$. If $n_{0,\alpha}=0$ we get   
$\alpha=-[i]$ 
for some $i\in\{1,\ldots,r\}$. 
 If $n_{0,\alpha}>0$, we apply $R(i)^*$ for an $i$ with $a_i\ne 0$. 
We see that all summands 
on the right side are divisible by $d_1=0$, and thus (P1) follows.

The relation
$R(m)$ is obtained from $R_{1,1,m,m}$ in two steps.
The relation $R_{1,1,m,m}$ immediately yields a
recursion relation identical to $R(m)$ except
for the fact that the sum is over $\vdash\!\!(d,\alpha)$ instead
of $\vdash\!\!(d,\alpha), d_i>0$.  
It will be shown that the terms with $d_1=0$ or $d_2=0$ vanish.
Since all summands are divisible by $d_1$, only the case $d_2=0$
need be considered. 
By (P1), either $N_{0,\gamma}=0$ or $\gamma=-[i]$.
In the second case, both binomal coefficients vanish.
Thus, relation $R(m)$ follows.

Now we show relation $R(i)$ holds.
We apply relation $R(i)^*$ to $N_{d,\alpha-[i]}$.
All summands on the right side of $R(i)^*$ are divisible by $d_1$,
thus all nonvanishing summands have $d_1>0$. By (P1),
$N_{0,\gamma}$ can only be nonzero if $\gamma=-[j]$ 
for some $j\in\{1,\ldots,r\}$. Since 
the right side of $R(i)^*$ is divisible by 
$c_i$,   the only nonzero summand on the right side with $d_2=0$
occurs for $(d_2,\gamma)=(0,-[i])$ and is $-d^2a_iN_{d,\alpha}$. 
Bringing this term on
the left side  and bringing 
$((a_i-1)^2-d^2)N_{d,\alpha-[i]}$ to the right side, 
we obtain the relation $R(i)$. Note
that $n_{d,\alpha}\ge 0$ implies   $n_{d,\alpha-[i]}\ge  1.$

We now show that 
the invariants $N_{d,\alpha}$ are determined recursively 
by the  relations  $R(1), \ldots, R(r)$, $R(m)$ and
the intial values. 
By (P1), all $d=0$ invariants are determined.
Let $d>0$.  
If $n_{d,\alpha}\ge 3$, then 
relation $R(m)$ determines $N_{d,\alpha}$ in terms of 
$N_{e,\lambda}$ with $e<d$.
Assume now that $0 \leq n_{d,\alpha}<3$. 
Either $(d,\alpha)=(1,(0, \ldots, 0))$ (and $N_{d,\alpha}=1$) or  
there exists an $i_0$ with $a_{i_0}\ne 0$.
By relation  $R(i_0)$, we can determine $N_{d,\alpha}$
in terms of $N_{e,\lambda}$ satisfying either $e<d$ or
$e=d$ and $n_{d,\lambda}>n_{d,\alpha}$. 
After at most 3 applications of a suitable $R(i)$, 
$R(m)$ may be
applied. $N_{d,\alpha}$ is then expressed in terms
of the intial values and  $N_{e, \lambda}$ with $e<d$.
This completes the recursion.

Finally, we verify   (P2)--(P5). First,
(P2) is proven. 
For $d=0$, the statement of (P2) is void. 
Let $d>0$, and assume by induction that 
(P2) holds for all $d_0<d$.  Let $(d,\alpha)$ be given with 
$d>0$, $a_j<0$. If  $n_{d,\alpha}\ge 3$, we can apply $R(m)$ to express
$N_{d,\alpha}$ as a linear combination of products 
$N_{d_1,\beta}N_{d-d_1,\alpha-\beta}$ with $d_1, d-d_1>0$. Furthermore 
$a_j<0$ implies $b_j<0$ or 
$a_j-b_j<0$. Therefore, $N_{d,\alpha}=0$ by induction.
If $0\leq n_{d,\alpha}<3$, we apply $R(j)$ to express 
$N_{d,\alpha}$ as a linear combination of $N_{d,\alpha-[j]}$ and
terms of the form $N_{d_1,\beta}N_{d-d_1,\alpha-[j]-\beta}$ with
$d_1, d-d_1>0$. These last terms vanish by induction. 
Thus $N_{d,\alpha}$ is just a multiple of $N_{d,\alpha-[j]}$.
As $n_{d,\alpha-[j]}=n_{d,\alpha}+1$, we can repeat this 
process to reduce to the case  $n_{d,\alpha}\ge 3$.

(P3) is obvious, as the initial values and  the set 
 $R(1), \ldots R(r), R(m)$ of relations are symmetric.

(P4) 
Let $(d,\alpha)$ be given. 
We will show that $N_{d,\alpha}=N_{d,(\alpha,0)}$.
By (P1) and the intial values, 
the result holds for $d=0$. Let $d>0$ and assume by induction
that the result holds for all $d_1<d$.
Case 1:  $n_{d,\alpha}\ge 3$. Apply $R(m)$ to express
$N_{d,\alpha}$ as a linear 
combination of terms  $N_{d_1,\beta}N_{d-d_1,\alpha-\beta}$ 
and to express
$N_{d,(\alpha,0)}$ 
as a linear combination of terms 
$N_{d_1,\beta_0}N_{d-d_1,(\alpha,0)-\beta_0}$
with $d_1, d-d_1>0$.
(P2) implies, for nonzero terms,
that $\beta_0$ must be of the form $(\beta,0)$.
Furthermore  the coefficient of 
$N_{d_1,(\beta,0)}N_{d_2,(\gamma,0)}$ in the expression for 
$N_{d,(\alpha,0)}$ is the same 
as that of $N_{d_1,\beta}N_{d_2,\gamma}$ in 
the expression for $N_{d,\alpha}$. 
Thus the result follows by induction on $d$.

Case 2: $0\leq n_{d,\alpha}<3$. 
If $\alpha\le 0$, then $(d,\alpha)$ must be $(1,(0\ldots,0))$ and
  $N_{d,\alpha}=N_{d,(\alpha,0)}=1$.
If there exists an $i$ with $a_i<0$, 
then $N_{d,\alpha}=N_{d,(\alpha,0)}=0$ by (P2). 
Assume there exists a $j$ with $a_j>0$.
We apply $R(j)$ both to $N_{d,\alpha}$ and $N_{d,(\alpha,0)}$.
$N_{d,\alpha}$ is expressed as a linear combination of 
$N_{d,\alpha-[j]}$ and the $N_{d_1,\beta}N_{d-d_1,\alpha-[j]-\beta}$ with 
$d_1, d-d_1>0$.  
Using (P2), the expression 
for  $N_{d,(\alpha,0)}$ is obtained by replacing  $N_{d_1,\beta}N_{d_2,\gamma}$
by  $N_{d_1,(\beta,0)}N_{d_2,(\gamma,0)}$ and $N_{d,\alpha-[i]}$
by $N_{d,(\alpha,0)-[i]}$. By induction on $d$,  it is enough 
to show the result for $N_{d,\alpha-[i]}$. Iterating the argument we  
reduce to $n_{d,\alpha}\ge 3$ or to  $\alpha\le 0$, 
where we already showed the
result.

(P5) Let $(d,\alpha)$ be given with 
$n_{d,\alpha}\ge 0$ and $a_j=1$ for some $j$. We show
that $N_{d, \alpha}=N_{d, \alpha -[j]}$. 
By (P1), we can assume $d>0$. 
We apply relation $R(j)$ to express $N_{d,\alpha}$ as a linear combination
of $N_{d,\alpha-[j]}$ and  
terms $N_{d_1,\beta}N_{d-d_1,\alpha-[j]-\beta}$ with 
$d_1, d-d_1>0$. Furthermore, by
(P2),  all nonzero terms have $b_j=c_j=0$. The coefficient of these terms is 
divisible by $c_j$. Therefore, $R(j)$ just reads
$d^2N_{d,\alpha}=d^2N_{d,\alpha-[j]}$.
\end{pf}

\section{Moduli Analysis}
\subsection{Results}
As before, let $X_r$ be the blow-up of $\P^2$
at $r$ general points $x_1, \ldots, x_r$.
In this section, the connection between  Gromov-Witten
invariants 
and the enumerative geometry of curves in $X_r$ is examined.
Let $\alpha=(a_1, \ldots, a_r)$. Let $(d,\alpha)$ denote
the class $dH- \sum_{i=1}^{r} a_i E_i$ in $H_2(X_r, \Z)$.
Let $n_{d,\alpha}= 3d-|\alpha|-1$ be the expected
dimension of the moduli space of maps $\overline{M}_{0,0}(X_r,
(d, \alpha))$.  If $n_{d,\alpha}\geq 0$, let $N_{d,\alpha}$
be the corresponding Gromov-Witten invariant. In this case,
the number of 
genus 0 stable maps of class $(d, \alpha)$ passing through
$n_{d,\alpha}$ general points of $X_r$ is proven to be
{\em finite}.
$N_{d,\alpha}$ is then
shown to be a count with (possible) multiplicities of the finite
solutions to this enumerative
problem.  
Hence, the Gromov-Witten invariant $N_{d,\alpha}$ is always
non-negative.
An analysis of the moduli space of maps yields a more
precise enumerative result.

\begin{thm} 
\label{numtheorem}
Let $n_{d,\alpha}\geq 0$, $d>0 $, and $\alpha\geq 0$.
Let (at least) one of the following two conditions hold for
the class $(d,\alpha)$: 
\begin{enumerate}
\item[(i)] $n_{d,\alpha} > 0$.
\item[(ii)] $a_i \in \{1,2\}$ for some $i$.
\end{enumerate}
Then,
$N_{d,\alpha}$ equals
the number of genus 0 stable maps of class $(d,\alpha)$ passing
through $n_{d,\alpha}$ general points in $X_r$. Moreover,
in this case, 
each solution map is an immersion of $\P^1$ in $X_r$.
\end{thm}

\subsection{Dimension 0 Moduli}
Three coarse moduli spaces of will be considered:
$$ M_{0,0}^\#(X_r, (d, \alpha)) \subset
M_{0,0}(X_r, (d,\alpha)) \subset \overline{M}_{0,0}(X_r, (d,\alpha)).$$
$M_{0,0}(X_r, (d,\alpha))$ is the open set of maps with domain
$\P^1$. $M_{0,0}^\#(X_r,(d,\alpha))$ is the open
set of maps with domain $\P^1$ that are {\em birational}
onto their image. As a first step, these unpointed
moduli spaces are shown to be empty when their expected dimensions
are negative. As always, $X_r$ is general.
 
\begin{lem} 
\label{lemone}
Let $(d,\alpha)\neq 0$ satisfy $n_{d,\alpha}<0$. Then, 
$\overline{M}_{0,0}(X_r, (d, \alpha))$
is empty.
\end{lem}
\begin{pf}
If $d<0$, $\overline{M}_{0,0}(X_r, (d, \alpha))$
is clearly empty.
Next, the case $d=0$ is considered. The only
classes $(0, \alpha)\neq 0$ that can be represented
by a connected curve are the classes $(0, -k[i])$ for $k\geq 1$.
Since $3\cdot 0 + k -1 \geq 0$, these
classes are ruled out by the assumption $n_{d,\alpha}<0$.
It can now be assumed that $d>0$.

Let ${\cal B}_r$ be the open configuration space
of $r$ distinct ordered points on $\P^2$. ${\cal B}_r$ is
an open set of 
$\P^2 \times  \cdots \times \P^2$ (with $r$ factors).
Let $\pi:{\cal X}_r \rightarrow {\cal B}_r$ be the universal family
of blown-up $\P^2$ 's. The fiber of $\pi$ over the
point $b=(b_1, \ldots, b_r)\in {\cal B}_r$ is simply $\P^2$ blown-up
at $b_1, \ldots, b_r$. The morphism $\pi$ is projective.
Let $\tau:\overline{{M}}_{0,0}(\pi, (d,\alpha))
\rightarrow {\cal B}_r$
be the relative coarse
moduli space of stable maps associated to
the family $\pi$. The morphism $\tau$ is projective. The
fiber $\tau^{-1}(b)$ is the corresponding moduli space of
maps $\overline{M}_{0,0}(\pi^{-1}(b), (d,\alpha))$ to the fiber 
$\pi^{-1}(b)$.

Assume that $\overline{M}_{0,0}(X_r, (d,\alpha))$ is nonempty
for general $X_r$. It follows that $\tau$ is a dominant
projective morphism and thus surjective onto
${\cal B}_r$. Let $b=(b_1,\ldots, b_r)\in {\cal B}_r$ be $r$ general
points on a nonsingular plane cubic $E \subset \P^2$.
Let $X_b=\pi^{-1}(b)$.
Since $\tau$ is surjective, there exists a stable
map $\mu: C \rightarrow X_b$.
By the numerical assumption, 
$$C \cdot \mu^*(c_1(T_{X_b}))= 3d-|\alpha| =n_{d,\alpha}+1\leq 0.$$
Since the points $b_1, \ldots, b_r$ lie on $E$,
the strict transform of $E$ is a representative of the
divisor class $c_1(T_{X_b})$ on $X_b$. Moreover,
since $E$ is elliptic, no component of $C$ surjects
upon $E$. Let $C=\bigcup C_j$ be the decompositon of
$C$ into irreducible components.
For each $C_j$, $\mu(C_j)$ is either a point
or an irreducible curve in $X_b$ not equal to $E$.
Hence, $C_j \cdot \mu^*(E) \geq 0$. Since
$$\sum_j C_j \cdot \mu^*(E) = C \cdot \mu^*(c_1(T_{X_b}))\leq 0,$$
$C_j \cdot \mu^{*}(E) =0$ for all components $C_j$.
Since $d>0$, there exists a component $C_{l}$
such that $\mu(C_{l})$ is of class $(d_{l}, 
\alpha_{l})$
with $d_{l}>0$.
Then, $\mu(C_{l})$ is curve and 
$\mu(C_{l}) \cap E= \emptyset$.
Now consider the image of $\mu(C_{l})$ in $\P^2$
(using the natural blow-down map $X_b \rightarrow \P^2$). The image of
$\mu(C_{l})$ is a degree $d_{l}>0$ plane curve meeting
$E$ only at the points $b_1, \ldots, b_r$. Hence, there is
an equality in the Picard group of $E$:
$${\cal O}_{\P^2}(d_l)|_E \stackrel{\sim}{=} 
{\cal O}_E(\sum_{i=1}^{r} m_i b_i)$$
for some non-negative integers $m_1, \ldots, m_r$.
Since $b_1, \ldots, b_r$ were chosen to be general
points on $E$, no such equality can hold. A contradiction
is reached and the Lemma is proven
\end{pf}

A map $\mu: \P^1 \rightarrow X_r$ is simply
incident to a point $y\in X_r$ if $\mu^{-1}(y)$
is scheme theoretically a single point in $\P^1$.

\begin{lem} 
\label{lemtwo}
Let $(d,\alpha)$ satisfy
$n_{d,\alpha}\geq 0$. 
Every map
 $[\mu]\in \overline{M}_{0,0}(X_r, (d,\alpha))$ incident
to $n_{d, \alpha}$ general points in $X_r$
is a birational map with domain $\P^1$. Morever, every
such map is simply incident to the $n_{d,\alpha}$ points.
\end{lem}
\begin{pf}
Let $C$ be a reducible curve.
Assume there exists a genus 0 (unpointed) stable map
$\mu:C \rightarrow X_r$
representing the class $(d, \alpha)$ incident
to $n_{d, \alpha}$ general points.
It is first claimed that at least two irreducible components are
mapped nontrivially by $\mu$. 
If no component is mapped to a point, the claim is trivial.
Otherwise, let $K$ be a maximal
connected component of $C$ that is mapped to a point.
$K$ must
meet the union of the irreducible components mapped nontrivially
in at least 3 points. Since $C$ is a tree, these
3 points lie on {\em distinct} components of $C$.
Let $C_1, \ldots, C_{s}$ be the irreducible
components mapped nontrivially by $\mu$.
Let $(d_1,\alpha_1), \ldots, (d_{s}, \alpha_{s})$ be the
classes represented by these components.
Let $p_i$ be the number of the $n_{d,\alpha}$ general
points contained in $\mu(C_i)$.
Since 
$$n_{d,\alpha}=s-1+\sum_{i=1}^{s}n_{d_i,\alpha_i}> 
\sum_{i=1}^{s}n_{d_i,\alpha_i},$$ 
and $\sum_{i=1}^{s} p_i \geq n_{d,\alpha}$, it follows that
for some $j$, $p_j>n_{d_j,\alpha_j}$. 
Let $y_1, \ldots ,y_{p_j}$ be the general points
contained in $\mu(C_j)$. Let $X_{r+p_j}$ be the
blow-up of $X_r$ at these points. Consider the strict
transform of the map $\mu$ to the map $\mu': C_j \rightarrow X_{r+p_j}$.
The class represented by $\mu'$ is $\beta=
(d_j, (\alpha_j, m_1, \ldots, m_{p_j}))$
where $m_i\geq 1$ for all $1\leq i \leq p_j$. Therefore
$n_\beta \leq n_{d_j, \alpha_j}-p_j <0$.
By Lemma (\ref{lemone}),
$\overline{M}_{0,0}(X_{r+p_j}, \beta)$ is empty.
A contradiction is reached.
Hence, no stable maps in $\overline{M}_{0,0}(X_r, (d,\alpha))$
with reducible domains pass through $n_{d,\alpha}$ general
points of $X_r$.

Next, assume there exists a stable map
$\mu: \P^1 \rightarrow X_r$ passing through $n_{d,\alpha}$
general points
which is not birational onto
its image. Let $\mu: \P^1 \rightarrow Im(\mu)$ be a generically
$k$-sheeted cover for $k\geq2$. 
Let $\gamma: \P^1 \rightarrow Im(\mu)$ be a
desingularization of the image.
The map $\gamma$ represents the class
$(d/k, \alpha/k) \neq (0,0)$ and is incident to the $n_{d,\alpha}$
general points.
Note that
$$n_{d/k, \alpha/k}=3\cdot {\frac{d}{k}} - {\frac{1}{k}} |\alpha| -1 <
n_{d, \alpha}.$$
As before, a contradiction is reached. 
Hence, the stable maps in $\overline{M}_{0,0}(X_r, (d,\alpha))$
passing through $n_{d,\alpha}$ general
points of $X_r$ are birational.

Finally, assume there exists a stable map
$\mu: \P^1 \rightarrow X_r$ passing through $n_{d,\alpha}$
general points $y_1, \ldots, y_{n_{d,\alpha}}$
which is not simply incident to the point $y_1$. 
Let $X_{r+n_{d, \alpha}}$ be the blow-up of $X_r$ at
the general points.
Then, the strict
transform of $\mu$ to $X_{r+n_{d,\alpha}}$ represents the class
$\beta=(d, (\alpha, m_1, \ldots, m_{n_{d,\alpha}}))$ where
$m_i\geq 1$ for all $1\leq i \leq n_{d,\alpha}$ and $m_1\geq 2$.
Again,
$n_\beta \leq n_{d,\alpha}- n_{d, \alpha}-1 <0$ and a contradiction is
reached.
\end{pf}

\begin{cor} 
\label{cortwo}
Let $(d,\alpha)$ satisfy
$n_{d,\alpha}=0$. Then $\overline{M}_{0,0}(X_r, (d,\alpha))
= {M}^\#_{0,0}(X_r, (d,\alpha))$.
\end{cor}

A scheme $Z$ is of {\em pure dimension 0} if
every irreducible component is a point. $Z$ may be empty. 

\begin{lem}
\label{lemthree}
Let $(d,\alpha)$ satisfy $n_{d,\alpha}=0$. Then, 
$\overline{M}_{0,0}(X_r,(d,\alpha))$ is of pure dimension 0.
\end{lem}
\begin{pf}
By Corollary (\ref{cortwo}), $\overline{M}_{0,0}(X_r, (d,\alpha)) =
M_{0,0}^\#(X_r, (d,\alpha))$.
Let $\mu: \P^1 \rightarrow X_r$ correspond to
a point $[\mu]\in M_{0,0}^\#(X_r, (d,\alpha))$.
Consider the normal (sheaf) sequence on $\P^1$
determined by $\mu$:
$$0 \rightarrow T_{\P^1} \rightarrow \mu^* 
T_{X_r} \rightarrow N_{X_r}
\rightarrow 0.$$ 
The sheaf $N_{X_r}$ has generic rank 1 and degree
equal to $3d-|\alpha|-2=n_{d,\alpha}-1=-1$.
There is a canonical torsion
sequence:
$$0 \rightarrow \tau \rightarrow N_{X_r} \rightarrow F \rightarrow
0.$$
The torsion subsheaf, $\tau$,
 is supported on the locus where
$\mu$ fails to be an immersion.
$F$ is a line bundle of degree equal to
$-1-dim(\tau)$.
It follows that 
\begin{equation}
\label{torr}
H^0(\P^1, N_{X_r})=
H^0(\P^1, \tau).
\end{equation}

Let $\lambda:{\cal C} \rightarrow M_{0,0}^\#(X_r, (d,\alpha))$
be any morphism of an irreducible curve to the moduli space. 
It will be shown that the
image of $\lambda$ is a point.
It can be assumed that ${\cal C}$ is nonsingular.
Since $M_{0,0}^\#(X_r, (d,\alpha))$ 
is contained in the automorphism-free locus, there exists
a universal curve $\pi: {\cal P} \rightarrow 
M_{0,0}^\#(X_r, (d,\alpha))$ and a universal
morphism $\mu: {\cal P} \rightarrow X_r$ (see [F-P]). Moreover,
$\pi$ is a $\P^1$-fibration. Let $\pi:S
\rightarrow {\cal C}$ be the pull-back of
${\cal P}$ via $\lambda$ and let $\mu:S \rightarrow X_r$
be the induced map. $S$ is a nonsingular surface.
Let $d\mu:T_S \rightarrow \mu^* T_{X_r}$ be the differential of $\mu$.
Let $T_V \subset T_S$ be the line bundle of $\pi$-vertical
tangent vectors, and let $U\subset S$ be the open
set where $d\mu: T_V \rightarrow T_{X_r}$ is a bundle
injection. The torsion result (\ref{torr}) directly
implies that the bundle map $d\mu: T_S \rightarrow T_{X_r}$
is of constant rank 1 on $U$. Hence, by the complex algebraic 
version of Sard's theorem, $\mu(S)$ is irreducible of
dimension $1$. 
The $\mu$-image of $S$ must equal
the $\mu$-image of each fiber of $\pi$. It now follows
easily that the image of $\lambda$ is a point.
\end{pf}

\subsection{The Map $\mu$ Over $E_i$}
The results of the previous section do not
show that $\overline{M}_{0,0}(X_r, (d,\alpha))$
is a nonsingular collection of points when $n_{d,\alpha}=0$.
Conditions for nonsingularity will be established in
section (\ref{nsns}).
Preliminary results concerning the the map $\mu$
over the exceptional divisors are required. First,
the injectivity of the differential over $E_i$ is
established.
\begin{lem}
\label{lemfour}
Let $(d,\alpha)$ satisfy $n_{d,\alpha}=0$.
Let $\mu: \P^1 \rightarrow X_r$ correspond
to a point $[\mu]\in \overline{M}_{0,0}(X_r, (d,\alpha))$.
Then $d\mu$ is injective at all points in $\mu^{-1} (E_i)$
(for all i).
\end{lem}
\begin{pf}
Consider again the relative coarse moduli space
$\tau:\overline{{M}}_{0,0}(\pi, (d,\alpha))
\rightarrow {\cal B}_r$ and the universal family
of blown-up $\P^2$'s,
$\pi:{\cal X}_r \rightarrow {\cal B}_r$.
Let ${\cal U}_r\subset {\cal B}_r$ denote the
open subset to which the conclusions of 
Corollary \ref{cortwo} and Lemma \ref{lemthree} apply.
For $b=(b_1, \ldots, b_r)\in {\cal B}_r$, 
let $E_i$ in $\pi^{-1}(b)$ denote
the exceptional divisor corresponding to the
point $b_i$.
Assume, for a general point $b\in {\cal U}_r$, there exists a
map $\mu: \P^1 \rightarrow \pi^{-1}(b)$ satisfying:
\begin{enumerate}
\item[(i)]
$[\mu] \in \overline{M}_{0,0}(\pi^{-1}(b), (d,\alpha))$.
\item[(ii)] There exists a point $p\in\P^1$ such that $d\mu(p)=0$
        and $\mu(p)\in E_i$ for some $i$.
\end{enumerate}
In this case,  there must
exist a fixed index $j$ such that for general $b\in {\cal U}_r$
the moduli space $\overline{M}_{0,0}(\pi^{-1}(b), (d,\alpha))$
contains a map with vanishing differential at some point over $E_j$.
Let $Y \subset \tau^{-1}({\cal U}_r)$ denote the locus
of maps with vanishing differential
at some point over  $E_j$. $Y$ is closed
in $\tau^{-1}({\cal U}_r)$. Let $\overline{Y}$ denote
the closure of $Y$ in 
$\overline{M}_{0,0}(\pi, (d,\alpha))$.
Let $[\mu] \in \overline{Y}$ where $\mu: C\rightarrow \pi^{-1}(\tau
([\mu]))$. It is 
easily seen that one of the following two cases hold:
\begin{enumerate}
\item[(i)] There exists a point $p\in C_{nonsing}$ satisfying
           $d\mu(p)=0$ and $\mu(p)\in E_j$.
\item[(ii)] There is a node of $C$ mapped to $E_j$.
\end{enumerate}
\noindent
These are the two possible degenerations of the singular point
of the morphism $\mu$ over $E_j$.
Since $Y$ dominates ${\cal B}_r$, the map $\overline{Y}
\rightarrow {\cal B}_r$ is surjective.

Define a complete curve ${\cal F}\subset {\cal B}_r$ as follows.
Let the points $e_1, \ldots, e_r$ be distinct points
on a nonsingular cubic plane curve $F \subset \P^2$. Choose
a zero for the group law on $F$.
Let the curve ${\cal F}\subset {\cal B}_r$ be determined
by elliptic translates of the tuple $(e_1, \ldots, e_r)$.
There is a natural map $\epsilon_j:{\cal F} \rightarrow F$
given by $\epsilon_j(f=(f_1,\ldots, f_r))=f_j$.
Consider the  fibration 
of blown-up $\P^2$'s over ${\cal F}$, 
$\pi^{-1}({\cal F})\rightarrow
{\cal F}$. Let $S\subset \pi^{-1}({\cal F})$ be the
subfibration of $\P^1$'s determined by the exceptional
divisor $E_j$.
$$S \subset \pi^{-1}({\cal F})  \rightarrow {\cal F}.$$
Via composition with $\epsilon_j$, there
is a natural projection $S\rightarrow F$. There is
a canonical isomorphism  $S \stackrel{\sim}{=} 
\P(T_{\P^2}|_F) \rightarrow F$ of varieties over $F$.

Let $\gamma:{\cal D} \rightarrow \overline{Y}$ be an irreducible
curve that surjects onto ${\cal F}$ via $\tau$.
After a possible
base change, a flat family of stable maps 
which induces the morphism $\gamma$ exists over ${\cal D}$. 
(In [F-P], the moduli space of maps is constructed locally
as finite quotient of a fine moduli space of
rigidified maps, so a base change with a universal
family exists on an open set of ${\cal D}$. The
properness of the functor of stable maps implies, after
further base changes, that this family can be completed over ${\cal D}$.)
Denote
this family of stable maps over ${\cal D}$ by
$\eta: {\cal C} \rightarrow {\cal D}$ and $\mu: {\cal C} \rightarrow
\pi^{-1}({\cal F})$. 
Let $Z \subset {\cal C}$ be the locus of 
nodes of the fibers of $\eta$ union the locus of
nonsingular points of the fibers where $d\mu$ vanishes
on the tangent space to the fiber. $Z$ is a closed subvariety.

Let $Z'\subset {\cal C}$ denote the (closed) intersection
$Z\cap \mu^{-1}(S)$. The subvariety $T=\mu(Z')\subset S= 
\P(T_{\P^2}|_F)$ dominates $F$ by the properties of $\overline{Y}$.
There is a natural section $F \rightarrow \P(T_{\P^2}|_F)$
given by the differential of $F$.
By Lemma \ref{lemfive}
below, $F\cap T$ is nonempty. Let $\zeta \in F \cap T$.

There are now two cases. First,
let  $d\in {\cal D}$ be  
such that there exists 
a nonsingular point $p\in{\cal C}_d$ at which the differential
of $\mu_d$ vanishes satisfing $\zeta=\mu_d(p)$. Consider the
map
$\mu_d$ from  ${\cal C}_d$ to $\P^2$ blown-up at the
points $f=(f_1, \ldots, f_r)$. 
Since $\zeta \in F \subset \P(T_{\P^2}|_F)$, 
the strict transform
of $F$ in this blow-up passes through $\zeta= \mu_d(p) \in E_j$.
If $p$ lies on a component of ${\cal C}_d$ not
mapped to a point, then
${\cal C}_d \cdot \mu^*(F) \geq 2$ because of the
vanishing differential at $p$.   
However, since $n_{d,\alpha}=0$ and $F$ represents the
first Chern class of the surface, ${\cal C}_d \cdot \mu^*(F) = 1$.
A contradiction is reached.
If $p$ lies on a component mapped to a point,
let $K$ be the maximal connected subcurve of ${\cal C}_d$ which contains
$p$ and is 
mapped to a point. By stability of the map, $K$ must intersect
the other components of ${\cal C}_d$ in at least 3 points.
By maximality, these intersection points lie
on components not mapped to a point by $\mu_d$. Hence, in this
case, ${\cal C}_d \cdot \mu^*(F) \geq 3$.    
Again a contradiction is reached.

Second,
let  $d\in {\cal D}$ be 
such that a node $p\in{\cal C}_d$ maps to $\zeta$. Again consider the
map
$\mu_d$ from  ${\cal C}_d$ to $\P^2$ blown up at the
points $f=(f_1, \ldots, f_r)$. The strict transform
of $F$ in this blow-up passes through $\zeta= \mu_d(p) \in E_j$.
If the node $p$ is an intersection of $2$ components
of ${\cal C}_d$ neither of which
is mapped to a point by $\mu_d$, then
${\cal C}_d \cdot \mu^*(F) \geq 2$ and a contradiction is reached.
If the node is on a component that is mapped to a point,
then ${\cal C}_d \cdot \mu^*(F) \geq 3$ as before and 
a contradiction is again reached.
\end{pf}

\begin{lem}\label{lemfive}
Let $\iota:F\hookrightarrow \P^2$ be a nonsingular plane cubic.
Let $F\rightarrow \P(T_{\P^2}|_F)$ be the canonical
section induced by the differential. Then $F \cap V$
is nonempty for any curve $V \subset \P(T_{\P^2}|_F)$.
\end{lem}
\begin{pf}
First the divisor class of the section $F$ is
calculated. Consider the tangent sequence
on the plane cubic $F$:
\begin{equation}
\label{fruit}
0 \rightarrow {\cal O}_F=T_F \rightarrow T_{\P^2}|_F \rightarrow 
{\cal O}_{\P^2}(3)|_F={\cal O}_F(3) \rightarrow 0.
\end{equation}
Let $S=\P(T_{\P^2}|_F)$ and let $\rho:S \rightarrow F$
denote the projection. Let $L$ denote the
line bundle ${\cal O}_{\P}(1)$ on $S$. Via
a degeneracy locus computation, sequence 
(\ref{fruit}) implies that the section $F$ is a divisor in
the linear series of the line bundle 
$L\otimes \rho^*{\cal O}_F(3)$. Note that:
$$H^0(S, L\otimes \rho^*{\cal O}_F(3))= H^0(F, T^*_{\P^2}|_F (3)).$$
The dual of the Euler sequence tensored with ${\cal O}_{\P^2}(3)$
restricted to $F$ yields:
$$0 \rightarrow T^*_{\P^2}|_F (3) \rightarrow 
\oplus_{1}^{3} {\cal O}_F(2) \rightarrow {\cal O}_F(3) \rightarrow 0.$$
It is easy to see the corresponding sequence on global
sections is exact. Hence $H^0(S, L\otimes \rho^*{\cal O}_F(3))=9$.
Therefore, for any $s\in S$, there exists a divisor
linearly equivalent to $F$ passing through $s$. Also,
it is easy to calculate $F\cdot F =9$.

Let $V$ be an irreducible curve in $S$ and assume $V\cap F$ is empty.
Hence, $V\cdot F=0$ and $V$ is not a fiber of $\rho$.
Let $G$ be a divisor equivalent to $F$ meeting $V$. By the
equation $V\cdot G=0$,
$V$ must be a component of $G$. Write $G=c_VV +\sum_i c_iW_i$.
Let $f$ be a general fiber of $\rho$. 
$$c_v V\cdot f + \sum_i c_iW_i \cdot f =G \cdot f=1.$$ 
$V \cdot f\geq 1$ since $V$ is
not a fiber. Therefore, $V\cdot f=1$, $c_V=1$,  and $W_i\cdot f=0$.
This implies each $W_i$ is a fiber. Then, 
$$9=F\cdot F=F\cdot G =  \sum_i F\cdot c_i W_i=\sum_i c_i.$$
$V$ is therefore a section of ${\cal O}_S(F) \otimes \rho^*N$
where $N$ is degree $-9$ line bundle on $F$. 
Again $H^0(S,   {\cal O}_S(F) \otimes \rho^*N)=
H^0(F, T^*_{\P^2}|_F \otimes {\cal O}_F(3) \otimes N)$.
The latter is seen to be zero by the dual Euler
sequence argument. No such $V$ exists.
\end{pf}

The Lemma (\ref{lemfour}) showed the branches of the
image curve $\mu(\P^1)$ are nonsingular at their
intersections with the $E_i$.
Next, it is shown that distinct branches of the image
curve do not intersect in the exceptional divisors.
\begin{lem}
\label{lemextra}
Let $(d,\alpha)$ satisfy $n_{d,\alpha}=0$.
Let $\mu: \P^1 \rightarrow X_r$ correspond
to a point $[\mu]\in \overline{M}_{0,0}(X_r, (d,\alpha))$.
Let $I$ be the image curve $\mu(\P^1)$. Then the
set $I\cap E_i$ is contained in the nonsingular locus
of $I$ (for all i).
\end{lem}
\begin{pf}
The proof of this Lemma exactly follows the proof of
Lemma (\ref{lemfour}). 
If the assertion is false, a quasi-projective 
subvariety $W\subset 
\overline{M}_{0,0}(\pi, (d, \alpha))$ can be found where
the image curve has distinct branches meeting in $E_j$
(for a fixed index $j$). The closure $\overline{W}$ of
$W$ then surjects upon ${\cal B}_r$. 
Let $\mu:C \rightarrow X_b$ be a limit map
$[\mu]\in 
\overline{W}$. At least one of the following properties
must be satisfied:
\begin{enumerate}
\item[(i)] Distincts point of $C$ are mapped by $\mu$
           to the same point of $E_j$.
\item[(ii)] There exists a point $p\in C_{nonsing}$ satisfying
           $d\mu(p)=0$ and $\mu(p)\in E_j$.
\item[(iii)] There is a node of $C$ mapped to $E_j$.
\end{enumerate}
The same curve ${\cal F}\subset {\cal B}_r$ is considered.
Let $\gamma:{\cal D} \rightarrow \overline{W}$ be an irreducible
curve that surjects onto ${\cal F}$ via $\tau$.
As before, a curve in $T\subset S=
\P(T_{\P^2}|_F)$ can be found representing the points
on $E_j$ where the singularities occur.
Using Lemma (\ref{lemfive}), $F\cap T$ is non-empty. It is then deduced
that stable maps exist satisfying $\mu^* c_1(T_{X_b}) \geq 2$
as before. A contradiction is reached.
\end{pf}

\subsection{Nonsingularity Conditions}
\label{nsns}
The main nonsingularity result needed
for the proof of Theorem (\ref{numtheorem}) can now be proven.
\begin{lem}
\label{lemmy}
Let $(d,\alpha)$ satisfy $d>0$, $\alpha\geq 0$, and
$n_{d,\alpha}=0$. 
If there exists an index $i$ for which $a_i \in \{1,2\}$,
then $\overline{M}_{0,0}(X_r, (d,\alpha))$
is nonsingular of pure dimension 0.
Moreover, the points of $\overline{M}_{0,0}(X_r, (d,\alpha))$
correspond to immersions of $\P^1$ in $X_r$.
\end{lem}
\begin{pf}
If $\overline{M}_{0,0}(X_r, (d,\alpha))$ is empty
for generic $X_r$, the Lemma is trivially true.
Let $\mu: \P^1 \rightarrow X_r$
be a map in $\overline{M}_{0,0}(X_r, (d,\alpha))$. By the
genericity assumption, the natural map:
\begin{equation}
\label{sirr}
d\tau: T_{ \overline{M}_{0,0}(\pi, (d,\alpha)), [\mu]}
\rightarrow  \tau^* T_{{\cal B}_r, \tau([\mu])}
\end{equation}
must be surjective. The 
Lemma is proved in two steps. First, the
surjectivity of (\ref{sirr}) is translated into
a condition on the global sections map of a normal
sheaf sequence associated to $\mu$.  
The map $\mu$ is then shown to be an {\em immersion}.
$N_{X_r}$ is therefore locally free of rank $1$ and
degree $3d-|\alpha|-2=n_{d,\alpha}-1 <0$.
The Zariski tangent space to  $\overline{M}_{0,0}(X_r, (d,\alpha))$
at $[\mu]$ is $H^0(\P^1, N_{X_r})=0$. Hence,
$[\mu]$ is a nonsingular point of 
$\overline{M}_{0,0}(X_r, (d,\alpha))$.

Let $X_r$ be the blow-up of $\P^2$ at
the points $x_1, \ldots, x_r$.
The deformation problem
as the blown-up points $x_1, \ldots, x_r$ vary is considered.
There is a projection $X_r \rightarrow \P^2$
which yields a sequence on $X_r$:
\begin{equation}
\label{www}
0 \rightarrow T_{X_r} \rightarrow
T_{\P^2} \rightarrow Q \rightarrow 0.
\end{equation}
$Q$ is a sheaf supported on the exceptional curves
$E_i$. $Q|_{E_i}$ is a line bundle on $E_i$.
More
precisely, if the point $e\in E_i$
corresponds to the tangent direction $T_e \subset T_{\P^2,x_i}$,
then the fiber of $Q$ at $e$ is $T_{\P^2,x_i}/ T_e$.
The space of deformations of the points $x_1, \ldots, x_r$
is $\oplus _{i=1}^{r} T_{\P^2, x_i} = H^0(X_r, Q)$.
$\oplus_{i=1}^r T_{\P^2, x_i}$ is also canonically
the tangent space to ${\cal B}_r$ at the
point $x=(x_1,\ldots, x_r)$. Therefore a vector
$0 \neq v\in \oplus_{i=1}^r T_{\P^2, x_i}$
defines a first order deformation of $X_r$ in the
family ${\cal X}_r$. Let 
$\lambda:\delt\rightarrow {\cal B}_r$ be a nonsingular
curve in ${\cal B}_r$ passing through $x$ with tangent direction 
$\C v$.
Let ${\cal X}_{\delt}= \lambda^{-1} {\cal X}_r \rightarrow \delt$.
This deformation naturally yields a differential sequence
on $X_r$:
\begin{equation}
\label{rrr}
0 \rightarrow T_{X_r} \rightarrow
T_{{\cal X}_{\delt}} \rightarrow {\cal O}_{X_r} \rightarrow 0.
\end{equation}
Sequences (\ref{www}) and (\ref{rrr}) are related
by a commutative diagram:
\begin{equation}
\label{sss}
\begin{CD}
0 @>>> T_{X_r} @>>> T_{{\cal X}_\delt} @>{a}>> {\cal O}_{X_r} @>>> 0 \\
@VVV   @VV{=}V  @VV{b}V @VV{c}V @VVV \\
0 @>>> T_{X_r} @>>>  T_{\P^2} @>{d}>>  Q @>>> 0.
\end{CD}
\end{equation}
Moreover, it is easy to check that the image of
$c:H^0(X_r,{\cal O}_{X_r}) \rightarrow H^0(X_r,Q)$ is
simply $\C v$.

Since $d\geq 1$, $Im(\mu)$ is not contained in
any $E_i$. Therefore the above commutatitive
diagram {\em stays exact} when 
pulled back to $\P^1$. Let $N_{\P^2}$ and 
$N_{{\cal X}_{\delt}}$ denote 
the normal sheaves on $\P^1$ of the maps to $\P^2$
and ${\cal X}_{\delt}$
induced by $\mu$. 
Consider the commutative diagram of exact sequences
obtained by pulling back (\ref{sss}) to $\P^1$ and
quotienting by the inclusion 
of sheaves induced by the differential $d\mu:T_{\P^1}\rightarrow
\mu^* T_{X_r}$.
\begin{equation*}
\begin{CD}
0 @>>> N_{X_r} @>>> N_{{\cal X}_\delt} @>{a}>> {\cal O}_{\P^1} @>>> 0 \\
@VVV   @VV{=}V  @VV{b}V @VV{c}V @VVV \\
0 @>>> N_{X_r} @>>>  N_{\P^2} @>{d}>>  \mu^*Q @>>> 0.
\end{CD}
\end{equation*}

$H^0(\P^1, N_{{\cal X}_\delt})$ is the space of first
order deformation of the map $\mu$ considered
as a map to ${\cal X}_{\delt}$. By the surjectivity
of (\ref{sirr}), there must exist a first
order deformation of $[\mu]$ not contained in $X_r$.
Therefore, the image
of $a: H^0(\P^1, N_{{\cal X}_{\delt}}) \rightarrow H^0(\P^1,
{\cal O}_{\P^1})$ must be non-zero. This condition
is equivalent to the splitting of the top sequence.
Using this splitting and the morphism $b$, it is seen that 
the section $v \in H^0(\P^1, Q)$ must be in the
image of 
$d: H^0(\P^1, N_{\P^2}) \rightarrow H^0(\P^1, \mu^*Q)$.

The conclusion of the above considerations is the
following. 
For every element $v\in \oplus _{i=1}^{r} T_{\P^2,x_i}$, there 
exists a section of $H^0(\P^1, N_{\P^2})$ which
has image $v\in H^0(\P^1, \mu^*Q)$.
The map $\mu$ will now be shown to be an immersion.

Suppose $p\in \P^1$ satisfies $\mu(p)\in E_i$.
By Lemma (\ref{lemfour}), $d\mu(p)$ is injective.
Let $m$ be the multiplicity of
$\mu^* E_i$ at $p$. Local calculations show that the
following hold in a neighborhood $U\subset \P^1$ of $p$
with local parameter $t$:
\begin{enumerate}
\item[(i)] $N_{\P^2}$ has torsion part $\C[t]/(t^{m-1})$
(where $t$ is a local parameter at $p$).
\item[(ii)] $\mu^*(Q)$ is the torsion sheaf $\C[t]/(t^m)$.
\item[(iii)] The map on torsion parts from $N_{\P^2}$ to
            $\mu^*(Q)$ is multiplication by $t$.
\end{enumerate}
Let $\tau$ be the torsion part of $N_{\P^2}$. By (iii), the natural
map of sheaves on $U$: 
$$N_{\P^2}/\tau \rightarrow \mu^*(Q) \otimes {\cal O}_p=\C$$
is surjective. Therefore,
a section $\overline{s}$ of the line bundle $N_{\P^2}/\tau$
is zero at $p$ if and only if the image of $\overline{s}$ in
$\mu^*(Q) \otimes {\cal O}_p$ is zero.

Decompose $\tau= A\oplus B$ where $A$ is the torsion
part supported at the points $\bigcup_i \mu^{-1}(E_i)$
and $B$ is the torsion part supported elsewhere.
Let $n$ equal the set theoretic cardinality
$|\bigcup_i \mu^{-1}(E_i)|$.
For each point $z\in \P^1$ lying over an exceptional
divisor $E$, let $m_z$ be the multiplicity of $\mu^* E$
at $z$.  The equations are obtained:
$$\sum_{z\in \bigcup_i \mu^{-1}(E_i)} m_z = \sum_i a_i,$$
$$degree(A)=
\sum_{z\in \bigcup_i \mu^{-1}(E_i)} (m_z-1) = -n+ \sum_i a_i.$$
The degree of $N_{\P^2}$ is $3d-2$. The degree of
$N_{\P^2}/A= 3d-2 +n -\sum_i a_i = n-1$.
Let $b=degree(B)$. Then, the degree of $N_{\P^2}/\tau$
is $n-1-b$. Note that $\mu$ is an immersion if and only
if $b=0$.

Without loss of generality, let $a_1\in\{1,2\}$.
First consider the case $a_1=1$. There is a unique
point $z_1$ in $\mu^{-1}(E_1)$. Let $v=\oplus_i v_i$
where $v_i\in T_{\P^2, x_i}$ 
satisfy:
\begin{enumerate}
\item[(i)]  $v_1\neq 0$ in $\mu^* Q \otimes {\cal O}_{z_1}$.
\item[(ii)] $v_i=0$ for $i\geq 2$. 
\end{enumerate}
Since there exists a section $s$ of $H^0(\P^1, N_{\P^2})$
with image $v\in H^0(\P^1, \mu^*(Q))$, there
must exist a nonzero section $\overline{s}$ of $H^0(\P^1, 
N_{\P^2}/\tau)$
vanishing at (at least) $n-1$ points (all the $z$'s except $z_1$)
by (iii). Therefore, $degree(N_{\P^2}/\tau) \geq n-1$. It follows
that $b=0$.

Next, consider the case $a_1=2$. There are two
possiblities. Either $\mu^{-1}(E_1)$ consists of
two points or one point. If there is a unique point
in $\mu^{-1}(E_1)$, the argument proceeds exactly
as in the $a_1=1$ case and $b=0$. Now suppose
$\mu^{-1} (E_1)= \{z_1, z_2\}$.
By Lemma (\ref{lemextra}), $\mu(z_1) \neq \mu(z_2)$.
 Let $v=\oplus_i v_i$ 
satisfy:
\begin{enumerate}
\item[(i)]  $v_1\neq 0$ in $\mu^* Q \otimes {\cal O}_{z_1}$.
\item[(ii)] $v_1=0$ in $\mu^* Q \otimes {\cal O}_{z_2}$.
\item[(iii)] $v_i=0$ for $i\geq 2$. 
\end{enumerate}
Such a selection of $v_1$ is possible since
$T_{\P^2, x_1}$ surjects upon $\mu^* Q \otimes {\cal O}_{z_1}
\oplus \mu^* Q \otimes {\cal O}_{z_2}$
for $\mu(z_1)\neq \mu(z_2)$.
As before, there
must exist a nonzero section $\overline{s}$ of $H^0(\P^1, N_{\P^2}/\tau)$
vanishing at least $n-1$ points (all the $z$'s except $z_1$)
by (iv). Therefore, $degree(N_{\P^2}/\tau) \geq n-1$. It follows
that $b=0$.
\end{pf}

\begin{lem}
\label{bee}
Let $d>0$, $\alpha\geq 0$, $r\leq 8$, and $n_{d, \alpha}=0$.
         Then, $\overline{M}_{0,0}(X_r,(d, \alpha))$ is
          nonsingular of pure dimension 0.
Moreover, the points of $\overline{M}_{0,0}(X_r, (d,\alpha))$
correspond to immersions of $\P^1$ in $X_r$.
\end{lem}
\begin{pf}
Let $\mu: \P^1 \rightarrow X_r$ be a map
in $\overline{M}_{0,0}(X_r,(d,\alpha))$.
       By Lemma \ref{lemfour}, $\mu$ is an immersion at the points
       of $\P^1$ mapping to the exceptional curves $E_i$.
       Suppose  $p\in \P^1$ is a point where $\mu$ is not an
       immersion ($\mu(p) \notin E_i$).  
       Since the number of blown-up points
       $x_1, \ldots , x_r$ is at most $8$, there is curve
       in the linear series $3H-\sum_{i=1}^{8} E_i$
       passing through $\mu(p)$.  
       Let $F$ denote this cubic (which may be
       reducible). There are now two cases.
       If $\mu(\P^1)$ is not contained in any component
               of $F$, then $\P^1 \cdot \mu^*(F) \geq2$ because
               $\mu$ is not an immersion at $p$. This is
               a contradiction since the numerical assumption
               implies $\P^1 \cdot \mu^*(F) =1$.
        If $\mu(\P^1)$ is contained in a component of $F$,
        then $d$ must equal 1,2, or 3 (since $\mu$ is
        birational).
        For these low degree cases, 
         $\overline{M}_{0,0}(X_r,(d,\alpha))$ is empty
             unless $a_i=1$ for some $i$. Then,
        Lemma \ref{lemmy} yields a contradiction.
        We conclude $\mu$ is an immersion and
        $\overline{M}_{0,0}(X_r, (d, \alpha))$ is nonsingular.
\end{pf}

\subsection{Proof of Theorem
(\ref{numtheorem})}
First, the case $n_{d, \alpha}=0$ is considered.
Since  $d>0$, $\alpha \geq  0$, and $a_i \in \{1, 2\}$
(for some $i$), 
Lemma \ref{lemmy} shows that 
$\overline{M}_{0,0}(X_r, (d,\alpha))$
is a nonsingular set of points. By Proposition \ref{fb}, $N_{d, \alpha}$
equals the number of points in
$\overline{M}_{0,0}(X_r, (d,\alpha))$.
Moreover, by Lemma \ref{lemmy},
the
points of  $\overline{M}_{0,0}(X_r, (d,\alpha))$
represent immersions of $\P^1$.
Theorem (\ref{numtheorem}) is established for classes
$(d,\alpha)$ satisfying $n_{d,\alpha}=0$.

Proceed now by induction on $n=n_{d, \alpha}$. 
If $n_{d,\alpha}>0$, consider the class
$(d, (\alpha,1))$ on $\P^2$ blown-up at $r+1$
points $x_1, \ldots, x_{r+1}$. 
Certainly, $n_{d,(\alpha,1)}=n -1$. 
By property (P5) of section \ref{qcb},
$$N_{d, \alpha} = N_{d, (\alpha,1)}.$$
The class $(d,(\alpha,1))$ satisfies condition (ii)
in the hypotheses of Theorem (\ref{numtheorem}).
By induction, $N_{d, (\alpha,1)}$ equals
the number of genus 0 stable maps of class
$(d, (\alpha,1))$ passing through $n_{d,\alpha}-1$
points $p_1, \ldots, p_{n-1}$ 
in $X_{r+1}$.  This is precisely equal
to the number of stable maps of class
$(d, \alpha)$ passing through  the $n_{d, \alpha}$
points $p_1, \ldots, p_{n-1}, x_{r+1}$ in $X_r$ by
Lemma \ref{lemtwo}. 
Since the solution curves are immersions in $X_{r+1}$, it
follows easily that the corresponding curves in
$X_r$ are also immersions. The proof of Theorem \ref{numtheorem}
is complete.

\section{Symmetries and Computations}
\subsection{The Cremona transformation}
\label{cremmy}
Let $p_1, p_2, p_3$ be 3 non-collinear points in $\P^2$.
Let $L_1, L_2, L_3$ be the 3 lines
determined by pairs of points where $p_i, p_j \in L_k$ for
distinct indices $i,j,k$.
Let $S$ be the blow-up of $\P^2$ at the points
$p_1, p_2,p_3$. Let $E_1,E_2,E_3$ be the
exceptional divisors of this blow-up.
Let $F_1, F_2, F_3$ be the strict transforms of the lines
$L_1, L_2,L_3$.
The $F_k$ 
are disjoint $(-1)$-curves on $S$ and can be blown-down.
The resulting surface is another projective plane
$\overline{\P}^2$. The blow-down maps are:
\begin{equation}
\label{crem}
\P^2 \stackrel{e}{\leftarrow} S 
\stackrel{f}{\rightarrow} \overline{\P}^2.
\end{equation}
This is the classical Cremona transformation
of the plane.
Let $q_1,q_2,q_3\in \overline{\P}^2$ be the points
$f(F_1),f(F_2), f(F_3)$. Let $H$ and $\overline{H}$
denote the hyperplane classes in $A_1(\P^2)$
and $A_1(\overline{\P}^2)$ respectively.
There are now 2 bases of $A_1(S)$
corresponding
to the two blow-downs:
$H, E_1, E_2, E_3$
and $\overline{H}, F_1, F_2, F_3$.
The relationship between these bases is:
\begin{align*}&dH-a_1E_1-a_2E_2-a_3E_3 =\\
&(2d-a_1-a_2-a_3)\overline{H}
- (d-a_2-a_3)F_1 - (d-a_1-a_3)F_2 -(d-a_1-a_2)F_3.
\end{align*}

Let $x_4, \ldots, x_r \in \P^2$ be additional
general points on $\P^2$ which correspond 
via the maps (\ref{crem}) to general points
$s_4,\ldots, s_r\in S$ and $y_4,\ldots, y_r\in \overline{\P}^2$.
The blow-up of $S$ at the points 
$s_4, \ldots, s_r$ may be viewed as a general blow-up 
of $\P^2$ at $p_1,p_2,p_3, x_4, \ldots, x_r$ 
or as a general blow-up of $\overline{\P}^2$
at $q_1, q_2, q_3, y_4, \ldots, y_r$. 
Let $G_4, \ldots, G_r$ denote the exceptional
divisors of the blow-up of $S$.

Since 
the class $dH- a_1 E_1 - a_2 E_2 -a_3 E_3 -\sum_{i=4}^r a_iG_i$
equals the class
$$(2d-a_1-a_2-a_3)\overline{H}
- (d-a_2-a_3)F_1 - (d-a_1-a_3)F_2 -(d-a_1-a_2)F_3-\sum_{i=4}^r a_iG_i,$$
the Gromov-Witten invariant
$N_{d, \alpha}$ on the blow-up of $\P^2$ equals
the invariant $N_{d', \alpha'}$ on the blow-up
of $\overline{\P}^2$ where
$$(d',\alpha')= (2d-a_1-a_2-a_3, 
(d-a_2-a_3,d-a_1-a_3,d-a_1-a_2,a_4,\ldots, a_r)).$$
It follows that $\overline{M}_{0,0}(X_r, (d, \alpha))$
is nonsingular if and only if $\overline{M}_{0,0}(X_r,
(d', \alpha'))$ is nonsingular. Therefore, 
$N_{d, \alpha}$ is enumerative if and only if 
$N_{d', \alpha'}$ is enumerative. 
The Cremona symmetry of the Gromov-Witten invariants
of $X_r$ is discussed in [C-M] from a slightly different
perspective.

For example, let $(d, \alpha)= (10, (4,4,3,3,3,3,3,3,3))=(10,(4^2,3^7))$
where the last equality is just notational convenience.
Then, $n_{10,(4^2, 3^7)}= 30 -29-1=0$. The class
$(10, (4^2,3^7)$ does not satisfying condition (i) or (ii)
of Theorem (\ref{numtheorem}).
Applying the Cremona transformation,
$(d', \alpha')= (9, (3,3,2,3^6))$.
Theorem (\ref{numtheorem}) applies
to $(d', \alpha')$. Therefore, the moduli space
$\overline{M}_{0,0}(X_r, (10, (4^2,3^7))$ is nonsingular
(and all points correspond to immersions).
$N_{10, (4^2,3^7)}=520$ is enumerative in this case.

\subsection{Tables} 
\label{tbls}
The arithmetic genus of the class
$(d, \alpha)$ on $X_r$ is determined by:
$$g_a(d, \alpha)= \frac{(d-1)(d-2)}{2} - \sum_{i=1}^{r}
\frac {a_i(a_i-1)}{2}.$$
The arithmetic genus of a reduced, irreducible curve
is non-negative. By Corollary \ref{cortwo},
$\overline{M}_{0,0}(X_r, (d, \alpha))$ is empty  
when $g_a(d,\alpha)<0$ and $n_{d,\alpha}=0$.
A simple reduction to the case of
expected dimension zero shows that
$N_{d, \alpha}=0$ if $g_a(d, \alpha)<0$.

If $a_i+a_j>d$ for indicies $i\neq j$, then
$N_{d, \alpha}=0$ unless $(d, \alpha)=(1,(1,1))$.
This follows again by a reduction to
the expected dimension zero case. Then,
Corollary \ref{cortwo} shows that   
$\overline{M}_{0,0}(X_r, (d, \alpha))$ is empty
(unless $(d, \alpha)=(1,(1,1))$)
by considering the intersection of
a map with the line in $\P^2$ connecting the
points $x_i$ and $x_j$.

In the first table below, Gromov-Witten
invariants $N_{d, \alpha}$ for $d\leq 5$ and $\alpha\geq 0$
are listed.
By properties (P3), (P4), and (P5), it suffices
to list the invariants for
ordered sequences $\alpha$ satisfying $\alpha\geq 2$.
Moreover, if $g_a(d,\alpha)<0$ or if
$a_i+a_j>d$, the invariant vanishes and is omitted
from the table. The invariants were computed by 
a Maple program via the recursive algorithm  of the proof 
of Theorem \ref{recurr}.

\begin{tabular}{|l|l|l|l|l|l|} \hline
$d=1$&2 &  3 & 4 & 5 & 5  \\
\hline 
$N_1=1$ & $N_2=1$ & $N_3=12$ & $N_4=620$ &$N_5=87304$ &$N_{5, (2^6)}=1$  \\
      &      & $N_{3,(2)}=1$ & $N_{4, (2)}=96$ & $N_{5, (2)}= 18132$ 
&$N_{5, (3)}=640$  \\
      &     &   & $N_{4, (2^2)}=12$  &  $N_{5, (2^2)}=3510$
&$N_{5, (3,2)}=96$ \\
      &   &  & $N_{4, (2^3)}=1$ & $N_{5, (2^3)}=620$ &$N_{5, (3,2^2)}=12$ \\
      &  &  &$N_{4,(3)}=1$  & $N_{5, (2^4)}=96$ &$N_{5, (3,2^3)}=1$  \\
      &  &  & & $N_{5, (2^5)}=12$ & $N_{5, (4)}=1$ \\
\hline
\end{tabular}
\vspace{+15pt}

\noindent
The Cremona transformation applied to the class $(5,(2,2,2))$
yields $N_{5,(2,2,2)}=N_{4,(1,1,1)}$. By Property
(P5), $N_{4,(1,1,1)}=N_{4}=620$.
The following table lists all the Gromov-Witten
invariants for degrees $6$ and $7$ which are
not obtained from lower degree numbers by the
Cremona transformation.
\vspace{+15pt}

\begin{tabular}{|l|l|l|} \hline
$d=6$ & 7 & 7\\ \hline
$N_6=26312976$ & $N_7= 14616808192$         & $N_{7,(3,2)}=90777600$\\
$N_{6,(2)}=6506400$ & $N_{7,(2)}=4059366000$& $N_{7,(3,2^2)}=23133696$ \\
$N_{6,(2^2)}=1558272$ & $N_{7,(2^2)}=1108152240$& $N_{7,(3,2^3)}= 5739856$ \\
$N_{6,(2^3)}=359640$ & $N_{7,(2^3)}=296849546$ &$N_{7,(3,2^4)}=1380648$\\
$N_{6,(2^4)}=79416$ & $N_{7,(2^4)}=77866800 $ &$N_{7,(3,2^5)}=  320160$  \\
$N_{6,(2^5)}=16608$ & $N_{7,(2^5)}= 19948176$&$N_{7,(3,2^6)}=71040$\\
$N_{6,(2^6)}=3240$ &  $N_{7,(2^6)}=4974460 $&$N_{7,(3,2^7)}=14928$\\
$N_{6,(2^7)}=576$ & $N_{7,(2^7)}= 1202355$&$N_{7,(3,2^8)}=2928$\\
$N_{6,(2^8)}=90$  & $N_{7,(2^8)}=280128$& $N_{7,(3^2)}=6508640$\\
$N_{6,(3)}=401172$ & $N_{7,(2^9)}=62450$& $N_{7,(4)}= 7492040$\\
$N_{6,(3,2)}=87544$ & $N_{7,(2^{10})}=13188$& $N_{7,(4,2)}= 1763415$\\
$N_{6,(4)}=3840$ & $N_{7,(3)}=347987200$& $N_{7,(5)}=21504$\\
\hline
\end{tabular}
\vspace{+15pt}

In [D-I], the Gromov-Witten invariants of $X_6$ are
computed. Our computation $N_{6,(2^6)}=3240$ disagrees with
[D-I]. We have checked our number using different
recursive strategies.

Let $(d, \alpha)$ be a class for which
all the hypotheses of Theorem \ref{numtheorem} and
Lemma \ref{bee} fail. Then, $r\geq 9$, $3d=|\alpha|+1$, and $\alpha \geq 3$.
Hence, $d\geq 10$. If $d=10$, then there are
only two possible values (up to reordering) for $\alpha$:
$(4^2, 3^7)$ or $(5, 3^8)$.
The invariant $N_{10, (4^2, 3^7)}$ was show to be enumerative
by the Cremona transformation in section \ref{cremmy}.
Applying the transformation to $(10, (5,3^8))$
yields $(9,(4,2^2,3^6))$. Hence, $N_{10,(5,3^8)}= N_{9,(4,2^2,3^6)}=90$
is enumerative by Theorem \ref{numtheorem}. We have shown all
invariants of degree $d\leq 10$ are enumerative.
The only invariants of
degree $11$ not proven to be enumerative by the methods of
this paper correspond to the classes
$(11,(5, 3^9))$ and $(11,(4^2,3^8))$. 
$N_{11,(5,3^9)}=707328$ and $N_{11,(4^2,3^8)}=2350228$.
It is not known to the authors whether non-trivial
multiplicities arise.

\end{document}